\documentclass[12pt]{article}
\usepackage[dvipdfmx]{graphicx}
\usepackage{amsthm,amsmath,amssymb}
\usepackage{latexsym}
\usepackage{graphics}
\usepackage{url}
\usepackage{placeins}
\usepackage{hyperref}
\usepackage{color}
\usepackage{bm}
\usepackage{comment}
\usepackage{statex}

\hypersetup{
     colorlinks=true,
     linkcolor=black,
     filecolor=blue,
     citecolor = black,
     urlcolor=cyan,
     }

\definecolor{blue}{rgb}{0,0,0.9}
\definecolor{red}{rgb}{0.9,0,0}
\definecolor{green}{rgb}{0,0.50,0.10}
\definecolor{violet}{rgb}{0.5804,0.0000,0.8275}  

\def\@themcountersep{}

\newtheorem{THEO}{Theorem}[section]
\newtheorem{ALGo}[THEO]{Algorithm}
\newtheorem{CONJ}[THEO]{Conjecture}
\newtheorem{COND}[THEO]{Condition}
\newtheorem{CORO}[THEO]{Corollary}
\newtheorem{DEFI}[THEO]{Definition}
\newtheorem{EXAMP}[THEO]{Example}
\newtheorem{FACT}[THEO]{Fact}
\newtheorem{HYPO}[THEO]{Hypothesis}
\newtheorem{LEMM}[THEO]{Lemma}
\newtheorem{PROB}[THEO]{Problem}
\newtheorem{PROP}[THEO]{Proposition}
\newtheorem{REMA}[THEO]{Remark}
\newcommand{\theo}{\begin{THEO}}
\newcommand{\algo}{\begin{ALGo} \rm}
\newcommand{\cond}{\begin{COND}}
\newcommand{\conj}{\begin{CONJ}}
\newcommand{\coro}{\begin{CORO}}
\newcommand{\defi}{\begin{DEFI} \rm}
\newcommand{\examp}{\begin{EXAMP} \rm}
\newcommand{\fact}{\begin{FACT}}
\newcommand{\hypo}{\begin{HYPO} \rm}
\newcommand{\lemm}{\begin{LEMM}}
\newcommand{\prob}{\begin{PROB} \rm}
\newcommand{\prop}{\begin{PROP}}
\newcommand{\rema}{\begin{REMA} \rm}
\newcommand{\etheo}{\end{THEO}}
\newcommand{\ealgo}{\end{ALGo}}
\newcommand{\econd}{\end{COND}}
\newcommand{\econj}{\end{CONJ}}
\newcommand{\ecoro}{\end{CORO}}
\newcommand{\edefi}{\end{DEFI}}
\newcommand{\eexamp}{\end{EXAMP}}
\newcommand{\efact}{\end{FACT}}
\newcommand{\ehypo}{\end{HYPO}}
\newcommand{\elemm}{\end{LEMM}}
\newcommand{\eprob}{\end{PROB}}
\newcommand{\eprop}{\end{PROP}}
\newcommand{\erema}{\end{REMA}}

\def\0{\mbox{\bf 0}}
\def\1{\mbox{\bf 1}}
\def\2{\mbox{\bf 2}}
\def\3{\mbox{\bf 3}}
\def\4{\mbox{\bf 4}}
\def\5{\mbox{\bf 5}}
\def\6{\mbox{\bf 6}}
\def\7{\mbox{\bf 7}}
\def\8{\mbox{\bf 8}}
\def\9{\mbox{\bf 9}}

\def\b{\mbox{\boldmath $b$}}

\def\e{\mbox{\boldmath $e$}}
\def\f{\mbox{\boldmath $f$}}

\def\s{\mbox{\boldmath $s$}}

\def\u{\mbox{\boldmath $u$}}
\def\v{\mbox{\boldmath $v$}}

\def\x{\mbox{\boldmath $x$}}
\def\y{\mbox{\boldmath $y$}}

\def\A{\mbox{\boldmath $A$}}
\def\B{\mbox{\boldmath $B$}}
\def\C{\mbox{\boldmath $C$}}

\def\F{\mbox{\boldmath $F$}}

\def\H{\mbox{\boldmath $H$}}

\def\P{\mbox{\boldmath $P$}}
\def\Q{\mbox{\boldmath $Q$}}

\def\W{\mbox{\boldmath $W$}}
\def\X{\mbox{\boldmath $X$}}
\def\Y{\mbox{\boldmath $Y$}}

\def\GC{\mbox{$\cal G$}}

\def\OC{\mbox{$\cal O$}}

\def\SC{\mbox{$\cal S$}}

\def\Real{\mbox{$\mathbb{R}$}}

\def\s0{\mbox{\scriptsize \boldmath $0$}}

\def\Real{\mathbb{R}}
\def\coneK{\mathbb{K}}

\def\coneN{\mathbb{N}}

\def\SymMat{\mathbb{S}}

\begin{document}

\title{ \large 
An Exceptionally Difficult Binary Quadratic Optimization Problem with Symmetry: a Challenge for The Largest Unsolved QAP Instance Tai256c\\
} 

%
%
%

\author{
\normalsize 
Koichi Fujii\thanks{NTT DATA Mathematical Systems Inc., Tokyo
 160-00016, Japan 
({\tt fujii@msi.co.jp}).}, \and \normalsize
Sunyoung Kim\thanks{Department of Mathematics, Ewha W. University, Seoul, 
	52 Ewhayeodae-gil, Sudaemoon-gu, Seoul 03760, Korea 
			({\tt skim@ewha.ac.kr}). 
The research was supported  by   NRF 2021-R1A2C1003810.
}, \and \normalsize
Masakazu Kojima\thanks{Department of Industrial and Systems Engineering,
	Chuo University, Tokyo 192-0393, Japan 
  	 ({\tt kojima@is.titech.ac.jp}).
}, \and \normalsize
Hans D. Mittelmann\thanks{School of Mathematical and Statistical Sciences, 
Arizona State University, Tempe, Arizona 85287-1804, U.S.A.
({\tt mittelma@asu.edu}).},  \and \normalsize
Yuji Shinano\thanks{Department of Applied Algorithmic Intelligence Methods (A$^2$IM), 
Zuse Institute Berlin, Takustrasse 7,14195 Berlin, Germany ({\tt shinano@zib.de}).
}
}

\date{\normalsize 
October 31, 2022, Revised: \today}
\maketitle







\begin{abstract}
\noindent
Tai256c is the largest unsolved quadratic assignment problem  (QAP) instance in 
QAPLIB. It is known that QAP tai256c can be converted into a 256 dimensional 
binary quadratic optimization problem (BQOP) with a single cardinality constraint 
which requires the sum of the binary variables to be 92. 
As the BQOP is much simpler than the original QAP, the conversion increases 
the possibility to solve the QAP.
Solving exactly the BQOP, however, is still very difficult. 
Indeed, a 1.48\% gap remains between the best known upper bound (UB) and lower bound (LB) 
of the unknown optimal value. This paper shows 
that the BQOP admits a nontrivial symmetry,
a property that makes the BQOP very hard  to solve.  Despite this difficulty, 
 it is imperative to decrease the gap in order to ultimately solve the BQOP exactly. 
To effectively improve the LB, we propose an efficient BB method that incorporates 
a doubly nonnegative relaxation,  
the orbit branching and the isomorphism pruning. 
With this  BB method, a new LB with 1.25\% gap is successfully obtained, 
and computing an LB with $1.0\%$ gap is shown to be still quite difficult.   
\end{abstract}






\section{Introduction} 

For a positive integer $n$, we let 
$N=\{1,\ldots,n\}$ represent a set of locations and also a set of facilities.
Given $n \times n$ symmetric matrices $\A =[a_{ik}]$ and 
$\B=[b_{j\ell}]$, 
the quadratic assignment problem (QAP) is 
stated as 
\begin{eqnarray}
\zeta^* = \min_{\pi} \sum_{i \in N}^n\sum_{k\in N}^n a_{ik}b_{\pi(i)\pi(k)}, 
\label{eq:QAP0}
\end{eqnarray}
where $a_{ik}$ denotes the flow between facilities $i$ and $k$, $b_{j\ell} = b_{\ell j}$ 
the distance between locations $j$ and $\ell$, and $(\pi(1),\ldots,\pi(n))$ a permutation of 
$1,\ldots,n$ such that 
$\pi(i) = j$ if facility $i$ is assigned to location $j$. 
We assume that 
the distance $b_{jj}$ from $j \in N$ to itself and the flow $a_{ii}$ from $i \in N$ to 
itself are both zero.

The QAP is NP-hard in theory, and solving exactly  large scale instances 
({\it e.g.}, $n \geq 40$) is very difficult in practice. 
To obtain an exact optimal solution,  we basically need two types of techniques. 
The first one is for computing 
heuristic solutions. 
Heuristic methods such as tabu search, genetic method
and simulated annealing have been developed for the QAP 
\cite{CONNOLLY1990,GAMB1997,SKORIN1990,TAILARD1991}. 
Those methods 
often attain  near-optimal solutions 
and occasionally  even  exact optimal solutions. 
The exactness is, however, not guaranteed in general.
The objective value $\bar{\zeta}$ obtained by those methods serves as 
an upper bound (UB) for the unknown optimal value $\zeta^*$. 
The second technique is 
to provide a lower bound (LB) $\underline{\zeta}$ for $\zeta^*$. If 
$\underline{\zeta}  = \bar{\zeta}$ holds, then we can conclude that 
$\underline{\zeta}  = \zeta^* = \bar{\zeta}$. Various relaxation methods 
\cite{ANSTREICHER2000,GILMORE1962,LAWLER1963,PRendl09,ZHAO1998} have been 
proposed for computing LBs. 
The two techniques mentioned above play 
essential tools in the branch 
and bound (BB) method for QAPs 
\cite{ANSTREICHER2002,CLAUSEN1997,GONCALVES2015,LAWLER1963,PARDALOS1997,ROUCAIROL1987}.


In this paper, we focus on the largest unsolved instance tai256c in 
QAPLIB \cite{BURKARD1997,QAPLIB}.  The main purpose of this paper is to investigate the challenge to solve 
the instance and provide an improved lower bound. Nissofolk et al. 
\cite{NISSOFOLK2016a,NISSOFOLK2016} showed that 
tai256c can be converted into a $256$-dimensional binary quadratic optimization 
problem (BQOP) with a single cardinality constraint $\sum_{i=1}^{256} x_i = 92$.
\begin{eqnarray}
\zeta^* & = & \min \left\{ \x^T\B\x : \x \in \{0,1\}^n \ \mbox{and } \sum_{i=1}^n x_i = 92 \right\}. 
\label{eq:BQOP0}
\end{eqnarray}
See Section 2.2. Here each feasible solution $\pi$ of QAP~\eqref{eq:QAP0}
is converted to a feasible solution $\x$ of BQOP~\eqref{eq:BQOP0}. 
They further transformed the BQOP~\eqref{eq:BQOP0} to a mixed integer convex 
quadratic program (MIQP) by  
the non-diagonal quadratic convex reformulation technique (NDQCR) developed 
in \cite{JI2012} for the quadratic knapsack problems. 
An LB = 44,095,032 (1.48$\%$ gap with respect to the best known UB  
44,759,294) was obtained by applying CPLEX to the resulting MIQP, where CPLEX terminated 
in almost 8.5 days since the node limit exceeded. 
In Section~\ref{section:NDQCR}, we reconfirm this difficulty, 
even when applying 
 the current version of CPLEX (IBM ILOG CPLEX Optimization Studio (version 22.1.1.0)~\cite{CPLEX}) 
and Gurobi Optimizer (version 11.0.0)~\cite{GUROBI}, which are expected to be significantly more powerful than 
the version of CPLEX used over a decade ago, to the MIQP. See also Remark~\ref{remark:NFQCR}. 
This demonstrates that the simple BQOP equivalent to tai256c  remains  difficult to solve. 

We show that BQOP~\eqref{eq:BQOP0} admits a nontrivial symmetry property, inherited 
from tai256c:\vspace{-2mm}
\begin{description}
\item[(a) ] The matrix $\B$ satisfies 
\begin{eqnarray}
\x_{\sigma}^T\B\x_{\sigma} = \x^T\B\x \ \mbox{for every } \x \in \{0,1\}^n 
\ \mbox{and } \sigma \in \GC,\label{eq:symmetry3} 
\end{eqnarray}
where $\GC$ denotes a subgroup of the symmetry group $\SC_n$ on $\{1,\ldots,n\}$ with $|\GC| = $2,048.\vspace{-2mm}   
\item[(b) ]
BQOP~\eqref{eq:BQOP0} has at least 1,024 distinct feasible solutions with the 
best known UB 44,759,294.
\end{description}

The size of BQOP~\eqref{eq:BQOP0}, 256, is not larger than 
quadratic 
unconstrained binary problem (QUBO) instances whose optimal 
solutions are known in the benchmark problem sets  \cite{MITTELMANN2023,BIQMAC}.
In fact, 
all QUBO instances with dimension less than 500 in the sets were solved exactly. BQOP~\eqref{eq:BQOP0} involves 
a single cardinality constraint 
$\sum_{i=1}^{256}x_i = 92$ in binary variables $x_i$ $(i=1,\ldots,256)$, 
which is expected to make solving 
BQOP~\eqref{eq:BQOP0} easier in comparison to QUBOs, 
since it considerably reduces the number of binary feasible solutions. 
Morevoer, it is straightforward to transform BQOP~\eqref{eq:BQOP0} into a QUBO 
by adding a penalty term 
$\lambda \big( \sum_{i=1}^{256} x_i - 92\big)^2$ to the objective quadratic form 
$\x^T\B\x$ with a sufficiently large $\lambda > 0$.  

Suppose that a BB method is applied to the BQOP with the best known UB. 
Then we (implicitly) construct an enumeration 
tree of its subproblems, where a subproblem is pruned  
whenever an LB of the subproblem 
not less than the best known UB of the BQOP (or an optimal solution of the subproblem) 
is obtained. In general, as the size of a subproblem involving a feasible 
solution with the best known UB becomes larger, the subproblem is harder to prune.  
As a result of the feature (b), the BB method cannot terminate in earlier stages since  
at least $1,024$ distinct feasible solutions with the 
best known UB are distributed over subproblems. 
Hence, 
 (b) primarily contributes to the difficulty of solving BQOP~\eqref{eq:BQOP0}.

To address the difficulty mentioned above in numerically solving BQOP~\eqref{eq:BQOP0}
and to compute a new LB better than the known ones,  this paper proposes \vspace{-2mm} 
\begin{description}
\item[(c) ] a BB method to show that the unknown optimal value $\zeta^*$ is not less than 
a given $\hat{\zeta}$. 
\vspace{-2mm} 
\end{description}
Here {\em a target} LB $\hat{\zeta}$ is chosen in the interval of the best known LB = 
$\underline{\zeta}$ = 44,095,032 and UB = $\bar{\zeta}$ = 44,759,294. We fix $\hat{\zeta}$ 
and $\bar{\zeta}$ 
before starting the BB method. The BB method terminates immediately after an LB not less 
than $\hat{\zeta}$ is obtained. 
The proposed BB method implements the Lagrangian doubly nonnegative (Lag-DNN)  
relaxation \cite{KIM2013,KIM2021} 
as a lower bounding procedure for subproblems.

Using this method, we compute a new LB 
44,200,000 ($1.25\%$ gap) in 39.2  days on a Mac Studio (20 CPU), 
and provide estimates on the amount of work (the number of 
subproblems to be solved and the execution time) for larger LBs. 
If we chose $\hat{\zeta}$ to be the best known UB $=\bar{\zeta}$, 
then $\bar{\zeta}$ would be proved to be the optimal value. 
In this case, $2.6\cdot10^{12}$ days would be required to 
solve $6.7\cdot10^{16}$ Lag-DNN relaxation  
subproblems of BQOP~\eqref{eq:BQOP0}. 
This is not an accurate estimate  and the execution time certainly  depends on a BB method including a lower bounding procedure,
 a branching rule and 
the computer used. Nevertheless,  it illustrates the extreme difficulty  of solving the BQOP. 

\subsection*{Contribution of the paper and existing results}

Our first contribution is to show and analyze the nontrivial 
symmetry property in BQOP~\eqref{eq:BQOP0} induced from tai256c. 
This BQOP is a   
simple, low-dimensional, and extremely difficult BQOP instance. 
As mentioned above, the nontrivial symmetry property~\eqref{eq:symmetry3} makes 
the BQOP 
 hard to solve.

The  second  contribution of this paper is a BB method to prove that the unknown 
optimal value $\zeta^*$ is not less than a given target LB $= \bar{\zeta}$. 
A BB method with a target LB was originally developed 
for large scale QAPs, which was successful to obtain improved lower 
bounds for some of the QAP instances in QAPLIB \cite{BURKARD1997,QAPLIB}
including sko100a,$\ldots, $ sko100f, 
tai80b, tai100b and tai150b (see \cite{MITTELMANN2023a}). 
The size of QAP tai256c, however, was too large to handle by the original BB method 
for the QAPs.
In the proposed method, we employ three effective techniques: The first one is 
the Lag-DNN relaxation \cite{KIM2013,KIM2021} subproblems of 
BQOP~\eqref{eq:BQOP0} for the lower bounding procedure. This relaxation is 
(almost) equivalent to a DNN relaxation \cite[Theorem 2.6]{ARIMA2018}, 
which is known as one of the 
strongest (numerically tractable) conic relaxations for combinatorial optimization problems 
\cite{ITO2017}. The second one is the 
standard orbital branching \cite{OSTROWSKI2011,PFETSCH2019} for reducing 
the size of the enumeration tree of 
subproblems to be solved. 
The third one is the isomorphism pruning \cite{MARGOT2002}. 
This technique works effectively to improve
the computational efficiency since the equivalence of some distinct subproblems 
occurs in the enumeration tree even after the orbital branching is applied. 
With this BB method, we computed an  LB with 1.46\% in 0.6 days on a Mac Studio (20 CPU), 
 generating 11,594 nodes. As mentioned, 
a slightly larger LB with 1.48\% gap is 
the best known one obtained by Nissofolk et al. \cite{NISSOFOLK2016} in 8.5 days 
when the node limit of CPLEX was reached. 
Furthermore,  we computed a new LB with 1.25\% in 39.2 days, generating 
1,077,353 nodes, and demonstrated
 the considerable difficulty in achieving an LB with  1.01\% gap. 
This can also be regarded as an important contribution. 

The orbital branching and the isomorphism pruning incorporated into our BB method can 
completely avoid applying 
the lower bounding and branching procedures to more than one isomorphic subproblem of BQOP~\eqref{eq:BQOP0}. 
In addition,  its size, $n=256$, is small, so we might expect to 
solve the BQOP relatively easily. However, in Section~\ref{section:difficulty}, we discuss another significant reason why 
solving BQOP~\eqref{eq:BQOP0} is extremely difficult compared to popular benchmark QUBOs with 
size $n \leq 256$. 

Exploiting symmetries 
of QAPs in their SDP relaxation was discussed in 
\cite{DEKLERK2010a,DEKLERK2012a,PERMENTER2020} (also \cite{BROSCH2022} 
in their DNN relaxation). However, those results are not relevant to 
the subsequent discussion of this paper. 


\subsection*{Outline of the paper}

In Section 2, 
we introduce key components that will be utilized in the subsequent sections,
including the conversion of QAP~\eqref{eq:QAP0} into BQOP~\eqref{eq:BQOP0} satisfying the symmetry property \eqref{eq:symmetry3}, the Lag-DNN  
relaxation and the Newton-bracketing (NB) method for solving the relaxation. 
In Section~\ref{section:Gurobi}, we present computational results using 
DABS (Diverse Adaptive Bulk Search, a genetic algorithm-based search algorithm)~\cite{NAKANO2023},
Gurobi 
and 
CPLEX.  
We show that state-of-the-art BQOP solver
Gurobi could not improve the known LB $\underline{\zeta}$,
demonstrating the difficulty of the problem.
In Section~\ref{section:newLowerBound}, we describe the BB method (c) in  detail  and 
report numerical results. 
We conclude the paper in Section \ref{section:concludingRemarks}.


\section{Preliminaries}

\label{section:prelim}

\subsection{Notation and symbols}

Let $N = \{1,\ldots,n\}$. We are mainly concerned with the BQOP~\eqref{eq:BQOP0} induced from tai256c.  In that 
case, $n = 256$. 
Let $\Real^n$ denotes the $n$-dimensional Euclidean space of column vectors $\x = (x_1,\ldots,x_n)$, 
and $\Real^n_+$ its nonnegative orthant $\{\x \in \Real^n : x_i \geq 0 \ (i \in N)\}$. 
For $\x \in \Real^n$, $\x^T$ is the transpose of $\x$. 
For each permutation $\sigma$ of $N$ and each $\x \in \{0,1\}^n \subset \Real^n_+$, 
$\x_{\sigma}$ denotes $\x' \in \{0,1\}^n$ such that $x'_j = x_{\sigma(i)}$ $(i \in N)$. 
$\Real^{m \times n}$ denotes the linear space of $m \times n$ matrices. 
$\SymMat^n$ denotes the linear space of $n \times n$ symmetric matrices with the 
inner product $\A \bullet \B = \sum_{i\in N}\sum_{j \in N}A_{ij}B_{ij}$ for $\A, \ \B \in \SymMat^n$, 
$\SymMat^n_+$ the convex cone of positive semidefinite matrices in $\SymMat^n$, and $\coneN^n$ the convex 
cone of matrices with nonnegative elements in $\SymMat^n$.  

\subsection{Conversion from tai256c to BQOP~\eqref{eq:BQOP0}} 

\label{subsection:conversion}

QAP~\eqref{eq:QAP0} can be rewritten with  an $n \times n$ matrix variable $\X$ as 
a quadratic optimization problem: 
\begin{eqnarray}
\zeta^* & = & \inf \left\{(\A\X\B) \bullet \X
:\X \in \Pi\right\},\label{eq:QAP1}  
\end{eqnarray}
where 
$\Pi$ denotes the set of $n \times n$ permutation matrix. 
We note that each feasible solution $\pi$ of QAP~\eqref{eq:QAP0}, which is a permutation 
of $N$, corresponds to an $\X \in \Pi$ such that 
$X_{ij} = 1$ iff $\pi(i) = j$ for every $(i,j) \in N \times N$. In case of tai256c, $n=256$ and 
$\A$ can be represented as 
$\A = \f\f^T$ for $\f = {\scriptsize \begin{pmatrix} \e \\ \0\end{pmatrix}}$, where $\e$ denotes the $92$-dimensional column vector of $1$'s. 
Hence, the objective function $(\A\X\B) \bullet \X$ of QAP~\eqref{eq:QAP1} can be rewritten  as 
\begin{eqnarray*}
(\A\X\B) \bullet \X = (\f\f^T\X\B) \bullet \X^T = (\X^T\f)^T\B(\X^T\f). 
\end{eqnarray*}
We then see that 
$\x = \X^T\f \in \{0,1\}^n \ \mbox{and } \sum_{i\in N} x_i = 92 \ \mbox{for every } \X \in \Pi$.   
Conversely, if $\x \in \{0,1\}^n \ \mbox{and } \sum_{i \in N} x_i = 92$, 
then $\x = \X^T\f$ for some $\X \in \Pi$. Therefore, QAP~\eqref{eq:QAP1} 
(hence QAP~\eqref{eq:QAP0}) is equivalent to  BQOP~\eqref{eq:BQOP0}. 

\rema The conversion from tai256c to BQOP~\eqref{eq:BQOP0} can also be carried out 
by the clone shrinking method in \cite[Theorem1]{FISCHETTI2012}. 
\erema

\subsection{Symmetry of the matrix $\B$} 

\label{subsection:symmetryB}

We computed  $\GC$ by a simple implicit enumeration of permutations $\sigma$ satisfying~\eqref{eq:symmetry3}, 
and found: \vspace{-2mm}
\begin{itemize}
\item $\left|\GC\right| = $2,048 \vspace{-2mm}
\item The best known feasible solution $\x^*$ of BQOP~\eqref{eq:BQOP0} 
with the objective value, which is equal 
 to the best known UB 44,759,294 for tai256c, 
 is further expanded 
to the set of feasible solutions $\{ (\x^*)_\sigma : \sigma \in \GC\}$ with the common 
objective value, where $\left|\{(\x^*)_\sigma : \sigma \in \GC\}\right| = 1024$; 
$(\x^*)_\sigma = (\x^*)_{\sigma'}$ can occur for distinct $\sigma \in \GC$ and $\sigma' \in \GC$.\vspace{-2mm}   
\end{itemize} 
Computing the group $\GC$ of permutations can also be carried out with
the software called Nauty \cite{MCKAY2010}. 
The symmetry of $\B$ is utilized in orbital branching (Section 4.2) and eliminating equivalent 
subproblems (Section 4.5) which are implemented in the BB method for solving BQOP~\eqref{eq:BQOP0}.   

\subsection{A Lagrangian doubly nonnegative (Lag-DNN) relaxation of a linearly constrained 
QOP in binary variables}

\label{eq:LagDNN}

We briefly present a Lag-DNN relaxation, 
which was originally proposed 
in \cite{KIM2013} combined with the the bisection-projection (BP) method for computing 
LBs of linearly constrained QOPs in binary variables. 
More recently, the  BP method was further 
enhanced to the Newton-bracketing (NB) method \cite{KIM2021} by replacing the bisection method with the Newton method for the largest zero 
of a continuously differentiable convex function $g : \Real \rightarrow [0,\infty)$ 
(see  Figure 1 and Section~\ref{eq:NewtonBracket})
In our proposed BB method, the NB method is used  
 for computing LBs of BQOP~\eqref{eq:BQOP0}. 
BQOP~\eqref{eq:BQOP0} as well as its 
subproblem $\mbox{BQOP}(I_0,I_1)$ presented in Section 4.1, 
are  special cases of a linearly constrained QOP in binary variables.
\begin{eqnarray}
\zeta & = & \inf\left\{\u^T\C\u : \u \in \{0,1\}^n, \ \F\u - \b s= \0, \ s = 1 \right\}, \label{eq:BQOP10} 
\end{eqnarray}
where $\C \in \SymMat^n$, $\F \in \Real^{m \times n}$ and $\b \in \Real^n$. 

 BQOP~\eqref{eq:BQOP10} is rewritten to strengthen the DNN relaxation by introducing slack variable 
vector $\v \in \{0,1\}^n$ for $\u \in \{0,1\}^n$:   
\begin{eqnarray}
\zeta & = & \inf\left\{\u^T\C\u : 
\begin{array}{l}
(\u,\v,s) \geq \0, \ (u_j + v_j -s )^2 = 0 \  (j \in N), \\[3pt] 
u_jv_j = 0 \  (j \in N), \
(\F\u - \b s)^T(\F\u - \b s)=0, \ s^2 = 1
\end{array}
 \right\} 
 \label{eq:BQOP11}
\end{eqnarray}
Introducing a penalty function (or a Lagrange function) 
\begin{eqnarray*}
L(\u,\v,s,\lambda) & = & \u^T\C\u + \lambda \big(\sum_{j\in N}(u_j+v_j-s)^2 
+ \sum_{j \in N}u_jv_j \\
& \ & 
+ (\F\u - \b s)^T(\F\u - \b s) \big) \ \mbox{for every } (\u,\v,s,\lambda) \geq \0,
\end{eqnarray*}
we consider a simple QOP 
\begin{eqnarray*}
\zeta(\lambda) & = & \inf\left\{L(\u,\v,s,\lambda) : (\u,\v,s) \geq \0, \  s^2 = 1\right\}, 
\end{eqnarray*} 
where $\lambda \geq 0$ denotes a penalty parameter (or a Lagrangian multiplier). We can prove that 
$\zeta(\lambda)$ converges to $\zeta$ as $\lambda \rightarrow \infty$. See \cite[Lemma 3]{KIM2013}. 
%
The sum $\u^T\C\u  + \sum_{j \in N}u_jv_j +  (\F\u - \b s)^T(\F\u - \b s) \big)$, which correspond the sum of the first, 
third and forth term of $L(\u,\v,s,\lambda)$ and form a quadratic form in $(\u,\v,s) \in \Real^{2n+1}$, can be represented 
as $\Q^1\bullet {\scriptsize \Bigg(\begin{pmatrix}\u \\ \v \\ s\end{pmatrix}\begin{pmatrix}\u \\ \v \\ s\end{pmatrix}^T\Bigg)}$ for some $\Q^1 \in \SymMat^{2n+1}$. 
Additionally, 
the second term $\lambda \big(\sum_{j\in N}(u_j+v_j-s)^2\big)$ can be written as 
$\lambda \Q^2\bullet {\scriptsize \Bigg(\begin{pmatrix}\u \\ \v \\ s\end{pmatrix}\begin{pmatrix}\u \\ \v \\ s\end{pmatrix}^T\Bigg)}$ for some $\Q^1 \in \SymMat^{2n+1}_+$ since it 
is a positive semidefinite quadratic form in $(\u,\v,s) \in \Real^{2n+1}$ 
for each fixed $\lambda \geq 0$ and 
linear in $\lambda \geq 0$ for each fixed $(\u,\v,s) \in \Real^{2n+1}$. Thus, letting 
$\Q_{\lambda} = (\Q^1 + \lambda \Q^2)$, we have 
\begin{eqnarray*}
L(\u,\v,s,\lambda) 
& = & 
\Q_\lambda 
\bullet \Bigg(\begin{pmatrix}\u \\ \v \\ s\end{pmatrix}\begin{pmatrix}\u \\ \v \\ s\end{pmatrix}^T\Bigg) \ \mbox{for every } (\u,\v,s,\lambda) \geq \0,  
\end{eqnarray*}
and
\begin{eqnarray*}
\zeta(\lambda) & = & \inf\left\{\Q_{\lambda} 
\bullet 
\Bigg( \begin{pmatrix}\u \\ \v \\ s\end{pmatrix}\begin{pmatrix}\u \\ \v \\ s\end{pmatrix}^T\Bigg)
 :
 \begin{array}{l} (\u,\v,s) \geq \0, \\
\H\bullet \Bigg(\begin{pmatrix}\u \\ \v \\ s\end{pmatrix}\begin{pmatrix}\u \\ \v \\ s\end{pmatrix}^T\Bigg)
= 1
\end{array}
\right\}, 
\end{eqnarray*} 
where $\H$ denotes the $(2n+1) \times (2n+1)$ matrix with elements $H_{ij} = 0$ 
$((i,j) \not=(2n+1,2n+1))$ and $H_{2n+1,2n+1} = 1$.   
By replacing ${\scriptsize \begin{pmatrix}\u \\ \v \\ s\end{pmatrix}\begin{pmatrix}\u \\ \v \\ s\end{pmatrix}^T}$ by 
a matrix variable $\W \in \SymMat^{2n+1}$, we obtain a Lag-DNN relaxation of BQOP~\eqref{eq:BQOP10}. 
\begin{eqnarray}
\eta(\lambda) & = & \inf\left\{ \Q_\lambda \bullet \W : 
\W \in \coneK, \ \ 
\H \bullet \W=1 \right\},  \label{eq:LDNN} 
\end{eqnarray}
where $\coneK = \SymMat^{2n+1}_+ \cap \coneN^{2n+1}$ denotes the $(2n+1)$-dimensional DNN 
cone. 
We note that a DNN relaxation of BQOP~\eqref{eq:BQOP11} can be written as 
\begin{eqnarray*}
\eta & = & \inf\left\{ \Q^1 \bullet \W : \W \in \coneK, 
\ \Q^2 \bullet \W = 0, \ \H \bullet \W = 1 \right\}. 
\end{eqnarray*}
The Lag-DNN relaxation~\eqref{eq:LDNN} is almost as strong as the DNN relaxation above in the sense that 
$\eta \geq \eta(\lambda)$ converges monotonically to $\eta$ as $\lambda \rightarrow \infty$. 
See \cite[Theorem 2.6]{ARIMA2018}. In Section 4, which presents
numerical results 
for the BB method applied to BQOP~\eqref{eq:BQOP0}, a value of
$\lambda = 10^8/\parallel \Q^1 \parallel$ is used.  

\subsection{The Newton-bracketing (NB) Method \cite{ARIMA2018,KIM2021}}

\label{eq:NewtonBracket}

Given $b_0 > \eta(\lambda)$, the NB Method applied 
to~\eqref{eq:LDNN} generates a sequence of intervals $[a_k,b_k]$ $(k=0,1,\ldots)$ which 
converges to $\eta(\lambda)$ monotonically. In this section, we briefly present how  
the sequence is generated. 
For more details, 
we refer to \cite[Section 3]{ARIMA2017}, \cite[Section 4]{ARIMA2018} and \cite{KIM2021}. 
Throughout this section, $\lambda > 0$ is fixed. 
The dual of DNN problem~\eqref{eq:LDNN} can be written as 
\begin{eqnarray}
y^* & = & \sup \left\{ y : \Q_\lambda - \H y = \Y \in \coneK^* \right\}, \label{eq:LDNNdual} 
\end{eqnarray}
where $\coneK^* = \SymMat^{2n+1}_+ + \coneN^{2n+1}$ (the dual of the $(2n+1)$-dim. 
DNN cone $\coneK = \SymMat^{2n+1}_+ \cap \coneN^{2n+1}$). 
By strong duality (see \cite[Lemma 2.3]{ARIMA2018}), $y^* = \eta(\lambda)$ holds. 
Define the function $g:\Real \rightarrow \Real$ as 
\begin{eqnarray*}
g(y) & = & \inf \{ \parallel \Q_\lambda - \H y - \Y \parallel : \Y \in \coneK^*\}\  \mbox{for every $y\in \Real$}, 
\end{eqnarray*}
which satisfies the following properties:\vspace{-2mm}
\begin{description}
\item{(i) } $g(y) \geq 0$ for every $y \in \Real$. \vspace{-2mm}
\item{(ii) } $\Y = \Q_\lambda - \H y \in \coneK^*$ ({\it i.e.}, $(y,\Y)$ is a feasible solution of~\eqref{eq:LDNNdual}) if and only if $g(y) = 0$. \vspace{-2mm}
\item{(iii) } $g(y^*) = 0$. \vspace{-2mm}
\item{(iv) } $g(y) = 0$ if $y \leq y^*$ (since $\H \in \coneK^*$). \vspace{-2mm}
\end{description}
\noindent
Hence, $y^*$ corresponds to the maximum zero of $g$.
Furthermore, $g: (y^*,\infty) \rightarrow \Real$ is convex and continuously differentiable (\cite[Lemma 4.1]{ARIMA2018}). 
Therefore, we can generate a sequence 
$(y_k)_{k=0}^\infty$  
converging monotonically to $\eta(\lambda)$ by applying the Newton method with a given initial 
point $y_0  > y^*$. See Figure 1. The function value $g(y)$ and the derivative $g'(y)$ 
at  $y= y_k > y^*$ is not given explicitly but can be computed by the accelerated proximal 
gradient (APG) method \cite{BECK2009}. 
This method also computes $(\Y^1_k,\Y^2_k) \in \SymMat^{2n+1} \times \coneN^{2n+1}$ 
which (approximately) satisfies $\Q_{\lambda} - \H y_k = \Y^1 + \Y^2$. 
and $\Y^1 \in \SymMat^{2n+1}_+$. 
We obtain $a_k \leq y^*$ by letting 
\begin{eqnarray*}
a_k & = & y_k + (n+1) \min\{0,\mbox{the minimum eigenvalue of $\Y^1_k$}\}.  
\end{eqnarray*}
See \cite[Lemma 3.1]{ARIMA2017}. In each iteration of the NB 
method, most of its execution time is consumed to evaluate 
$g(y_k)$, $g'(y_k)$, 
$\Y^1_k$ and $\Y^2_k$ by the APG method.

\begin{figure}
\begin{center}
\includegraphics[width=0.35\textwidth]{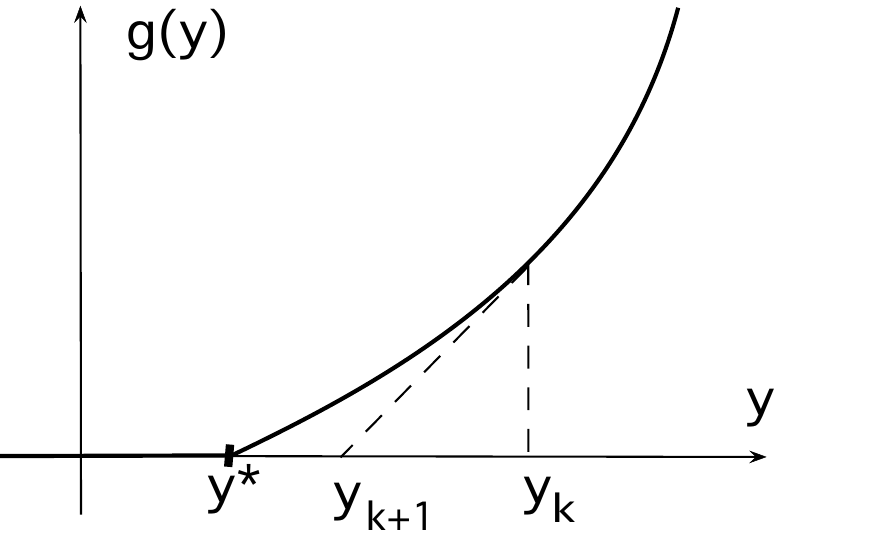} 
\caption{The convex function $g:\Real \rightarrow \Real_+$ and the NB method. 
}
\end{center}
\end{figure}


\section{Numerical experiments using the QUBO solver and MIQP solvers}
\label{section:Gurobi} 

\subsection{Solving BQOP~\eqref{eq:BQOP0} with QUBO and MIQP Solvers}

We report computational results obtained by 
  DABS (Diverse Adaptive Bulk Search), a genetic algorithm-based search algorithm for QUBO \cite{NAKANO2023},
and a general BQOP solver Gurobi Optimizer (version 11.0.0) \cite{GUROBI}, for 
 BQOP~\eqref{eq:BQOP0}. 
Numerical experiments were conducted on
Intel(R) Xeon(R) Gold 6246 (3.30 GHz) processors using 48 threads with 1.5TB of RAM.


  BQOP~\eqref{eq:BQOP0} with a single cardinality 
  constraint $\sum_{i=1}^n x_i = 92$ can be transformed to a simple QUBO 
  by adding the penalty term $\lambda(\sum_{i=1}^n x_i - 92)^2$ 
  to the objective function $\x^T\Q\x$ and removing the constraint, 
  where $\lambda > 0$ is a penalty parameter. 
We applied DABS  to the resulting QUBO with $\lambda=10^7$.
DABS attained a feasible solution 
  with the objective value of 44,759,294,
  which coincides with the best known upper bound \cite{BURKARD1997, QAPLIB}
  for tai256c, within a few seconds.
Moreover, the feasible solution computed is contained in the 1,024 known ones
  with the same objective value of 44,759,294 (see Section 2.3).

We also applied Gurobi to BQOP~\eqref{eq:BQOP0}. 
Gurobi is state-of-the-art as a solver for general BQOPs.
In a benchmark project conducted by Mittelmann~\cite{mittelmannbenchmarks}, 
several solvers are tested for BQOPs.
The results there indicate that Gurobi is the fastest
in solving those BQOPs.
Gurobi has been enhanced as a BQOP solver,
and its performance and efficiency in this specific area have been constantly
improved.
In Section~\ref{section:branching},
  we will see from the symmetry of $\B$ 
  that  BQOP \eqref{eq:BQOP0} has an optimal solution $\x$ with $x_1 = 1$. 
When  Gurobi was applied to BQOP \eqref{eq:BQOP0} with $x_1$ fixed to 1, the internal log indicated that 255 variables were grouped into 44 orbits, with one of the orbits consisting of only a single variable.
As we will see in Table\ref{table:orbits}(Section \ref{section:newLowerBound}), this coincides with
  the symmetry detected by our method.

To experiment with Gurobi, some parameters were needed to be decided, in particular 
\texttt{MIPFocus}, and \texttt{PreQLinearize}. 
The parameter \texttt{MIPFocus} controls the solution strategy of branch-and-bound. 
We chose \texttt{MIPFocus=3} which focuses on computing the LB. 
The parameter 
\texttt{PreQLinearize}
   controls presolving for BQOPs. More precisely,  
the parameter \texttt{PreQLinearize=0} 
adds neither any variables nor  constraints, 
  but it performs adjustments on quadratic objective functions
  to make them positive semidefinite.  
The parameter values  \texttt{PreQLinearize=1} and \texttt{PreQLinearize=2}
  attempt to linearize quadratic constraints or a quadratic objective,
  replacing quadratic terms with linear terms
  with additional variables and linear constraints.

In the first step, we examined the value of \texttt{PreQLinearize}
by comparing the LB 
obtained at the root node.
\texttt{PreQLinearize=0} provided the 
  LB 41,172,797,
  \texttt{PreQLinearize=1}
    the LB 10,759,778, and 
   \texttt{PreQLinearize=2} the LB 
   3,987,504.  
As the best LB  
41,172,797 among the three
was obtained by  \texttt{PreQLinearize=0},
 we adopted this setting for our experiments.

In the second step, we examined 1,024 potentially best solutions
generated by the symmetry of BQOP~\eqref{eq:BQOP0}. 
We determine the best initial solution among them
  by executing branch-and-bound with a time limit of 1 hour
  and comparing the LBs.  

In the final step, we executed branch-and-bound with 
\texttt{MIPFocus=3}, \texttt{PreQLinearize=0} and 
the initial solution chosen in the second step.
  We had to set a time limit of 60,000 seconds
  due to limitations in computational resources.
  We finally obtained an LB of  41,669,052
  using Gurobi, which generated 12,518,148 nodes. 

The LB 41,669,052 
obtained is significantly lower than the LBs we will present in Section~4. This suggests that, despite its rapid improvement in solving BQOP problems, the state-of-the-art general BQOP solver is still considerably less effective in solving 
BQOP~\eqref{eq:BQOP0} or improving its LB.  

There are two primary factors 
for Gurobi's poor performance in solving BQOP~\eqref{eq:BQOP0}.  
Recall that the LB 41,172,797 (8.01\% gap) was obtained at the root node problem. 
We applied the Lag-DNN relaxation to the same problem 
and computed an LB 43,881,304 (1.96\% gap) of the problem by the NB method 
in less than 10 minutes. This shows that 
the LB procedure incorporated in Gurobi is much weaker than the Lag-DNN relaxation. 
Another noteworthy factor contributing to Gurobi's poor performance
is its inability to fully utilize the 
symmetry property of BQOP~\eqref{eq:BQOP0}. 
This also played a role in the generation of a huge number of nodes, reaching 
12,518,148, by Gurobi to compute the LB 
41,669,052 (6.90\%) 
which is still much smaller than the LB 43,881,304 (1.96\%) 
of the root node computed by the Lag-DNN relaxation
%
This sharply contrasts with our improved BB method, 
which generated 1,077,353 nodes to compute the LB 44,200,000 (1.25\%) in 39.2 days. See Section 4.5. 

\subsection{NDQCR method} 

\label{section:NDQCR}

We also tested the NDQCR method \cite{NISSOFOLK2016a, NISSOFOLK2016} with Gurobi and IBM ILOG CPLEX Optimization Studio (version 22.1.1.0).
Numerical experiments were conducted on Intel(R) Xeon(R) Gold 6230 (2.10 GHz) processors using 32 threads with 385 GB of RAM.
The NDQCR involves 
the following three steps to construct an MIQP: \vspace{-2mm}
\begin{description}
\item{1.} 
Add a set of quadratic inequalities to BQOP~\eqref{eq:BQOP0} 
by applying the reformulation-linearization technique (RLT) \cite{SHERALI2013}. 
\vspace{-2mm}
\item{2.} Solve the primal-dual pair of SDP relaxations of the resulting quadratically constrained 
quadratic program (QCQP). 
\vspace{-2mm}
\item{3.} Formulate an MIQP effectively utilizing the information from the optimal solution of the 
dual SDP relaxation problem. 
\vspace{-2mm}
\end{description}
We formulated the dual problem of the SDP relaxation of BQOP(\ref{eq:BQOP0}), 
and solved it using MATLAB R2023b in 212.53 seconds.
We then solved the MIQP, which was strengthened by utilizing the optimal solution of the dual SDP relaxation, with CPLEX and Gurobi.
For both solvers, we conducted runs for 24 hours and 9 days, utilizing 32 threads.
With CPLEX, the LB 43,928,939 was obtained after 24 hours and the LB 43,950,303 after 9 days, reaching a total of 62,563,872 nodes.
With Gurobi, 
an LB  44,011,626 was obtained after 24 hours and an LB 44,035,000 after 4.97 days, exhausting 385 GB of memory and reaching a total of 139,859,181 nodes.
The NDQCR method can obtain the LB 43,848,767(2.03\% gap) at the root node, which is competitive with the NB method. However, the LB increases only gradually, causing the branch-and-bound tree to expand rapidly. 
The NDQCR method  
is expected to require increasingly longer CPU time and encounter difficulties with the rapid explosion of nodes needed to compute higher quality LBs. 

\rema 
\label{remark:NFQCR} 
The technical details of the NDQCR method tested above are different from the original 
NDQCR method proposed in \cite{NISSOFOLK2016a, NISSOFOLK2016} because the code of 
the original NDQCR method was not available from the authors. We tried to incorporate the basic 
ideas 1, 2 and 3 mentioned above to simulate the original NDQCR method. 
Although the NDQCR method tested did not attain the 
same performance as the original one, it encountered the same 
difficulties with the rapid explosion of nodes. 
\erema


\section{A branch and bound method for a given target lower bound}

\label{section:newLowerBound}

The optimal value $\zeta^*$ of BQOP~\eqref{eq:BQOP0} remains unknown, with
an LB $\underline{\zeta} =$ 44,095,032 $\leq \zeta^* \leq$ a UB $\bar{\zeta} =$ 44,759,294 exhibiting a $1.48\%$  
gap between them. 
We need to improve  the UB and/or the LB to compute the optimal value $\zeta^*$.
To improve the LB,  we propose a BB  method.  
For the lower bounding procedure, we use the Lag-DNN relaxation of a subproblem 
and the NB method for computing
its optimal value which serves as an LB of the subproblem. 
See Sections 2.4 and 2.5 for the Lag-DNN relaxation and the NB method, respectively. 
Any upper bounding procedure is not incorporated. 
Before the start of the BB method, a target LB, $\hat{\zeta}$, is first set such that 
$\underline{\zeta} = 44,095,032 < \hat{\zeta} \leq \bar{\zeta} = 44,759,294$. 
A target LB, $\hat{\zeta}$, is the desired value to obtain.  
Ideally, we want to set $\hat{\zeta} = \bar{\zeta}$ to confirm whether  $\bar{\zeta}$ 
is the optimal value. But such a setting may be too ambitious, 
and requires much stronger computing power than
the machine currently used.  
As a larger $\hat{\zeta}$ is set,
the computational cost rapidly increases as we will see in Section~\ref{section:estimation}. 

We describe a class of 
subproblems of BQOP~\eqref{eq:BQOP0} which appear in the enumeration tree generated by the BB method in 
Section~\ref{section:subproblms}, and the orbital branching, the branching procedure used in the BB method in Section~\ref{section:branching}. 
Before presenting numerical results on the BB method in Section~\ref{section:result}, we 
provide  a preliminary estimate in Table 2 for the amount of work (the number of nodes to generate and execution time) 
to attain given target LBs by the BB method in Section~\ref{section:estimation}. 
We note that the preliminary estimate was obtained 
 using the orbital branching, without employing the isomorphism pruning  
 described in Section 4.5. 
Based on this estimation, 
we choose appropriate LBs as targets 
for the numerical experiment on the BB method whose results 
are reported in Section~\ref{section:result}. In Section~\ref{section:isomorphism}, we present 
 the isomorphism pruning \cite{MARGOT2002} 
to improve the performance of the BB method. 
Our numerical results on the original and improved BB methods, which are 
summarized in Tables 3 and 4, respectively, demonstrate that 
employing the isomorphism pruning 
is expected to 
halve the computational effort compared to using |red{the orbital} branching alone.
A new LB $\hat{\zeta}=$ 44,200,000 (1.25\% gap) 
is also attained by the improved BB method in 39.2 days. 

\subsection{A class of subproblems of BQOP~\eqref{eq:BQOP0}} 

\label{section:subproblms}

Let 
\begin{eqnarray*}
\SC & = & \left\{ (I_0,I_1,F) : 
\begin{array}{l}
\mbox{a partition of $N$, {\it i.e.},} \
I_0 \bigcup I_1\bigcup F = N, \\ 
I_0, \ I_1 \ \mbox{and } F \ \mbox{are disjoint with each other}
\end{array}
\right\}. 
\end{eqnarray*}
Obviously, $F$ is uniquely determined by $I_0$ and $I_1$ as 
$F = N \backslash \big(I_0 \bigcup I_1 \big)$ for each $(I_0,I_1,F) \in \SC$. Hence, 
$F$ in the triplet $(I_0,I_1,F)$ is redundant, and we frequently omit $F$ for the simplicity 
of notation. 
For each $(I_0,I_1,F) \in \SC$, 
we consider a subproblem of BQOP~\eqref{eq:BQOP0} 
\begin{eqnarray*}
\mbox{BQOP}(I_0,I_1) \ : \ \zeta(I_0,I_1) & = & \min
\left\{ \x^T \B \x : 
\begin{array}{l}
\x \in \{0,1\}^n, \ \displaystyle \sum_{i=1}^n x_i = 92, \\
x_i = 0 \ (i \in I_0), \ x_j = 1 \ (j \in I_1)
\end{array} 
\right\}\\
& = & \min \left\{ \y^T \B(I_0,I_1) \y : 
\begin{array}{l}
\y \in \{0,1\}^{F}, \\
 \displaystyle \sum_{i \in F} y_i = 92 - \left| I_1 \right|
\end{array} 
\right\} 
+\sum_{j\in I_1}\sum_{k\in I_1} B_{jk},
\end{eqnarray*}
where 
\begin{eqnarray*}
& & \y \in \Real^F \ \mbox{denotes the subvector of $\x$ with elements $x_i$ $(i\in F)$}, \\ 
& & \B(I_0,I_1) = \B_{FF} + 
2 \times \mbox{diagonal matrix of} \left(\sum_{k\in I_1} \B_{kF} \right), \\ 
& &
\B_{kF} = \mbox{the column vector consisting of elements $B_{kj}$ $(j \in F)$
}.
\end{eqnarray*}
To see the equivalent representation of $\B(I_0,I_1)$ above, 
we note that each $x_k$ $(k \in F)$ in BQOP($I_0,I_1)$  is a binary variable; hence $x_jx_k= x_k^2$ if 
$j \in I_1$ and $k \in F$. 
For example, if $F = \{1,\ldots,\ell\}$ and $I_1\cup I_0 = \{\ell+1,\ldots,n\}$, 
then $\B(I_0,I_1)$ is an $\ell \times \ell$ matrix with elements 
$B(I_0,I_1)_{ij}$ $(i=1,\ldots,\ell,j=1,\ldots,\ell)$ such that 
\begin{eqnarray*}
B(I_0,I_1)_{ij} = 
\left\{
\begin{array}{ll}
B_{ij} & \mbox{if $i \not= j$}, \\ 
B_{ii} + 2 \sum_{k \in I_1} B_{ki} & \mbox{if $i=j$}. 
\end{array}
\right.
\end{eqnarray*}
For computing an LB of $\mbox{BQOP}(I_0,I_1)$ in the BB method,  
we applied the NB method to the Lag-DNN relaxation of $\mbox{BQOP}(I_0,I_1)$ with 
$\lambda = 10^8 / \parallel \B(I_0,I_1) \parallel$. 

\subsection{Orbital branching} 

\label{section:branching}

We discuss the orbital branching technique 
from~\cite{OSTROWSKI2011,PFETSCH2019}.
As  mentioned in Section 1, $\B = \B(\emptyset,\emptyset)$ satisfies the 
symmetry property~\eqref{eq:symmetry3}. This property is partially shared  by   
numerous $\B(I_0,I_1)$ instances where $((I_0,I_1,F) \in \SC)$. Let $(I_0,I_1,F) \in \SC$ be fixed.
Assume in general that 
\begin{eqnarray}
\y_{\sigma}^T \B(I_0,I_1) \y_{\sigma} = \y^T \B(I_0,I_1) \y \ 
\mbox{for every } \y \in \{0,1\}^{|F|} \ \mbox{and } \ \sigma \in \GC(I_0,I_1)  
\label{eq:symmetry52}
\end{eqnarray}
holds, where $\GC(I_0,I_1)$ is a group of permutations of $F$. 
Let 
$\omega(i) = \{j \in F: j = \sigma_i \ \mbox{for some } \sigma \in \GC(I_0,I_1)\}$ for every $i \in F$, and $\OC(I_0,I_1) = 
\left\{\omega(i): i \in F \right\}$. 
Each $o \in \OC(I_0,I_1)$ is called an {\em orbit} of the  
group $\GC(I_0,I_1)$. 
Let $\min(o)$ denote the minimum index of orbit~$o$, which 
serves as a representative for~$o$. 
Then 
we know that all BQOP$(I_0,I_1\bigcup \{j\})$ 
$(j \in o)$ are equivalent in the sense that they share a common optimal value 
$\zeta(I_0,I_1\bigcup \{\min(o)\})$. 
Therefore, we can branch BQOP$(I_0,I_1)$ to two sub BQOPs, 
BQOP$(I_0\bigcup o,I_1)$ and BQOP$(I_0,I_1\bigcup $ 
$ \{\min(o)\})$. 

In general, $\OC(I_0,I_1)$ consists of multiple orbits. Selecting an appropriate  $o$ 
from $\OC(I_0,I_1)$ for branching of BQOP$(I_0,I_1)$ to 
BQOP$(I_0\bigcup o,I_1)$ and BQOP$(I_0,I_1\bigcup \{\min(o)\})$
is an important issue to design an efficient branch and bound method. In our 
numerical experiment presented in Section~\ref{section:result}, 
an orbit $o$ is chosen from $\OC(I_0,I_1)$ 
according to the average objective value of 
BQOP$(I_0,$ $I_1\bigcup\min(o))$  over all 
feasible solutions, so that the chosen orbit, $o^*$,  attains the largest value.  
Then, we apply the branching of BQOP$(I_0,I_1)$ to 
two subproblems BQOP$(I_0$ $\bigcup o^*,I_1)$ and BQOP$(I_0,I_1\bigcup \{\min(o^*)\})$. 
Here the average objective value 
is computed as 
$\x^T\B\x$ with 
$x_i = 0 \ (i\in I_0), \ x_j = 1 \ (j \in I_1)$ and $x_k = (92-|I_1|) / |F| \ (k \in F)$. 

$\GC(\emptyset,\emptyset) = \GC$  has the single orbit 
$o = N = \{1,\ldots,n\}$. We branch BQOP$(\emptyset,\emptyset)$ 
into two subproblems BQOP$(N,\emptyset)$ and BQOP$(\emptyset,\{1\})$. 
Obviously, the former BQOP$(N,\emptyset)$ is infeasible. 
Table 1 
summarizes the branching of the node BQOP$(\emptyset,\{1\})$ to  
BQOP$(\{2,16,17,$ $241\},\{1\})$ and BQOP$(\emptyset,\{1,2\})$, where orbit $\{2,16,17,241\}$ 
is chosen from $\OC(\emptyset,\{1\})$. 

\begin{table}[htp]
\scriptsize{
\begin{center}
\label{table:orbits}
\begin{tabular}{|c|l|c|c|c|c|}
\hline
Orbit  &       & The size & The average objective \\
number & \multicolumn{1}{c|}{Orbit } & of orbit & value of BQOP$(\emptyset,\{1,\min(o)\})$ \\
\hline
 1   & 2 16 17 241    &  4  & 52655297.0 \\ 
 2   & 18 32 242 256    &  4  & 52567852.0 \\
 3   & 3 15 33 225    &  4  & 52524130.0 \\
 4   & 19 31 34 48 226 240 243 255    &  8  & 52515385.0 \\
 5   & 35 47 227 239    &  4  & 52502268.0 \\
     & $\cdots$ & & $\cdots$ \\ 
30   & 87 91 102 108 166 172 183 187    &  8  & 52483274.0 \\
31   & 9 129    &  2  & 52483139.0 \\
32   & 72 74 117 125 149 157 200 202    &  8  & 52483097.0 \\
33   & 25 130 144 249    &  4  & 52483097.0 \\
    & $\cdots$ & & $\cdots$ \\  
43   & 121 136 138 153    &  4  & 52481955.0 \\
44   & 137    &  1  & 52481773.0 \\
\hline
\end{tabular}
 \vspace{3mm}
 \caption{
A summary of branching of BQOP$(I_0,I_1)$ 
with $I_0 = \emptyset$, $I_1 = \{1\}$ and $F = \{2,3,\ldots,256\}$ 
to BQOP$(\{2,16,17,241\},\{1\})$ and BQOP$(\emptyset,\{1,2\})$. 
Here $F = \{2,3,\ldots,256\}$ is partitioned into 44 orbits, 
which consist of $21$ orbits with size $8$, 
$21$ orbits with size $4$, $1$ orbit with size $2$ and $1$ orbit with size $1$.  
The $44$ orbits are listed according to the decreasing order of the average 
objective value of BQOP$(\emptyset,\{1,\min(o)\})$ over all feasible solutions.
}
\end{center}
}
\end{table}

In addition to the branching rule mentioned above, we employ 
the simple breadth first search; the method to search the enumeration tree is not relevant 
to the computational efficiency since the incumbent objective value is fixed to 
the target LB $\hat{\zeta}$ and any upper bounding procedure is not applied. 
At the initial ($0$th) iteration, we take BQOP($I_0,I_1$) with $I_0 = \emptyset$ and $I_1 = \{1\}$ 
as the active root node.
Suppose that 
at the start of the $k$th iteration with $k \geq 0$, all active nodes located at  depth $k$ of the enumeration tree have already  been generated.
The $k$th iteration consists of two phases.  First, we compute LBs of 
all active nodes to determine whether they remain active by applying the NB method. Second, we apply 
the orbital branching to the resulting active nodes, which are then located at depth $(k+1)$ of the enumeration tree, 
for the $(k+1)$th iteration.
At each node BQOP$(I_0,I_1)$ of the enumeration tree, 
the NB method generates a sequence of intervals $[a_p,b_p]$ $(p=1,2,\ldots)$ 
satisfying a monotonicity property: 
(1) $a_p$  converges monotonically to an LB $\nu$ of BQOP$(I_0,I_1)$
from below, and (2) $b_p$  converges monotonically to $\nu$ from above. 
Thus, if $\hat{\zeta} \leq a_q$ holds for some $q$, we know that the LB $\nu$ to which 
the interval $[a_p,b_p]$ converges is not smaller than $\hat{\zeta}$, 
and BQOP$(I_0,I_1)$ can be pruned. 
On the other hand, if $b_q < \hat{\zeta}$ holds for some $q$, we know the LB $\nu$ 
of BQOP$(I_0,I_1)$ is 
smaller than $\hat{\zeta}$; hence  the iteration can be stopped and  branching  
to BQOP$(I_0,I_1)$ can be applied. Therefore, the above properties (1) and (2) 
of the NB method 
work very effectively to increase the computational efficiency of the BB method. See Figure 3. 

\subsection{Estimating the total number of nodes generated by the BB method}

\label{section:estimation}

All the computations for numerical results reported in this section and the next two sections were 
performed using MATLAB 2022a on a Mac Studio with Apple M1 Ultra CPU, 20 cores 
and 128 GB memory.
For the parallel computation, we solved Lag-DNN relaxations of 20 subproblems BQOP$(I_0,I_1)$ 
in parallel by the NB method with  the `parfor' loop of MATLAB.  

\begin{table}[htp]
\scriptsize{
\begin{center}
\begin{tabular}{|r|r||r|r|r|r|r|r|r|r|r|}
\hline
\multicolumn{2}{|c||}{target LB}    &  \multicolumn{3}{c|}{no. of nodes} &  \multicolumn{3}{c|}{exec. timed (day)} \\
$\hat{\zeta}$ & gap & min & mean & max & min & mean & max \\
\hline
44,100,000 & 1.46\% & 22,175 & 27,278 & 29,692 &  0.9 &  1.1 &  1.2  \\ 
\hline
44,120,000 & 1.43\%  & 53,625 & 72,574 & 91,944 &  2.1 &  2.8 &  3.6  \\ 
\hline
44,130,000  & 1.41\% & 64,084 & 133,417 & 275,264 &  2.5 &  5.2 &  10.7  \\ 
\hline
44,150,000  & 1.36\%  & 241,827 & 293,696 & 339,245 &  9.4 & 11.5 &  13.2  \\ 
\hline
44,200,000  & 1.25\%  & 827,791 & 1,983,516 & 2,891,498 & 32.3 & 77.4 &  112.8  \\ 
\hline
44,300,000  & 1.03\% & 5.4$\cdot10^7$ & 3.7$\cdot10^8$ & 8.1$\cdot10^8$ & 2.1$\cdot10^3$ & 1.4$\cdot10^4$  & 3.1$\cdot10^4$  
\\ 
\hline
44,500,000  & 0.58\% & 1.3$\cdot10^{11}$ & 5.5$\cdot10^{12}$ & 2.6$\cdot10^{13} $ & 5.1$\cdot10^{6}$ & 2.2$\cdot10^8$ & 1.0$\cdot10^9$  \\ 
\hline
44,759,294 & 0.00\%   & 1.2$\cdot10^{14}$ & 6.7$\cdot10^{16}$ & 3.3$\cdot10^{17}$ & 4.6$\cdot10^{9}$ & 2.6$\cdot10^{12}$ 
&  1.3$\cdot10^{13}$  \\ 
\hline
\end{tabular}
\vspace{3mm}
\caption{Estimation of the work of the BB method described in Sections 4.1 and 4.2. 
For each target LB $\hat{\zeta}$, we applied 
5 different random sampling of $s_k$ nodes from $\bar{t}_k$ active nodes at 
the depth $k$ with $k \geq \ell$ for the next depth $k+1$. Here min, mean and max 
denote the minimum, the mean 
and the maximum of  those 5 estimations of the number of nodes in the enumeration tree 
to generate and the execution time (day), respectively. 
A Lag-DNN subproblem at each node was solved in about 30 $\sim$ 150 seconds.} 
\end{center}
}
\end{table}

To choose a reasonable target LB  \ $\hat{\zeta}$ which can be  attained by the BB method, 
we performed 
preliminary numerical experiments to estimate the computational work. 
Given a target LB $\hat{\zeta}$, 
we construct an enumeration tree by the breadth first search as long as the number $t_k$ 
of nodes at the depth $k$ of the tree is smaller than $1000$. 
Suppose that $t_0, t_1,\ldots,t_{\ell-1} <  1000 \leq t_{\ell}$; 
hence the full enumeration tree has been constructed up to the depth $\ell$ 
by the BB method. We start sampling at the depth $\ell$ and construct a random subtree to estimate the total number 
of nodes in the full enumeration tree.
Let $\bar{t}_{\ell} = t_{\ell}$. 
At each depth $k \geq \ell$, we choose $s_k$ nodes randomly from $\bar{t}_k$ active nodes 
for the next depth $(k+1)$, where 
\begin{eqnarray*}
	s_k & = & \left\{
	\begin{array}{ll}
	100 & \mbox{if } \bar{t}_{k} \geq 500, \\
	\bar{t}_{k} & \mbox{otherwise}.  
	\end{array}
\right. 
\end{eqnarray*}
Then, we apply the lower bounding procedure using the NB method to the selected $s_k$ 
nodes and the branching procedure to the resulting $r_k$ active nodes to generate a subset 
of the nodes in the full enumeration tree at the depth $(k+1)$. 
Next, we let $\bar{t}_{k+1} = 2  r_k$, 
which is the cardinality of the subset (the number of nodes in the subset) 
as each active node is branched into two child nodes. 
We may regard $2r_k /s_k = \bar{t}_{k+1}/s_k$ 
as the increasing rate of the nodes from the depth $k$ to the depth ${k+1}$, and the total 
number of nodes in the full enumeration tree is estimated by 
\begin{eqnarray}
& & \sum_{k=1}^{\ell} t_k  + \sum_{k > \ell} \hat{t}_{k}, \ \mbox{where } \
    \hat{t}_{\ell} = t_{\ell}, \ \hat{t}_{k+1} = (2r_k/s_k)\hat{t}_k \ (k \geq \ell). 
\label{eq:noOfNodes}
\end{eqnarray}
We continue this process till $r_k$ attains $0$.  
Table 2 shows the estimation of computational work (the number of nodes to generate and 
the execution time) for $8$ cases $\hat{\zeta}=$ 44,100,00,$\ldots$,44,759,294.
In spite of the simplicity of this unrefined method, 
it provides useful information on whether a given target LB 
can be attained by the BB method on the computer used. 

\subsection{Numerical results}

\label{section:result}

We see from Table 2 that the cases with the target LB 
$\hat{\zeta} = 44,759,294$, $44,500,000$ and $44,300,000$ are very challenging. 
The case $\hat{\zeta} = 44,200,000$ could be processed but might take more than a 
few months. 
Table 3 shows numerical results for the other 4 cases with $\hat{\zeta}$ = 44,100,000, 44,120,000, 44,130,000 and 44,150,000. 
We observe that the estimation of the number of nodes and execution time described in the previous section 
are useful. 

\begin{table}[htp]
\scriptsize{
\begin{center}
\begin{tabular}{|r|r||rl|rl|r|}
\hline
\multicolumn{2}{|c||}{target LB} & \multicolumn{2}{c|}{no. of nodes} &  \multicolumn{2}{c|}{exec. time (day)} \\
$\hat{\zeta}$ & gap &  & estimation (min,mean,max) & &  estimation (min,mean,max)\\ 
\hline
44,100,000  & 1.46\% & 23,510 & (22,175, 27,278, 29,692) &  1.0 & ( 0.9, 1.1, 1.2) \\
44,120,000  & 1.43\% & 63,554 & (53,625, 72,574, 91,944) &  2.5 & ( 2.1, 2.5, 3.6) \\
44,130,000  & 1.41\% & 102,310 & (64,084, 133,417, 275,264) &  4.1 & ( 2.5, 5.2, 10.7) \\
44,150,000  & 1.36\% & 277,304 & (241,827, 293,696, 339, 245) &  10.7 & (9.4, 11.5, 13.2) \\
\hline
\end{tabular}
\vspace{3mm}
\caption{Numerical results on the BB method described in Sections 4.1 and 4.2. The 3 numbers in the parenthesis ($\cdot,\cdot,\cdot$) 
denote the minimum,  
 mean 
 and maximum estimation from Table 2, respectively. 
} 
\end{center}
}
\end{table}

Figure 2 (A) displays 
the change in the number of nodes  
as the depth $k$ increases, and 
Figure 2 (B)  the change in 
the number of nodes with size 2 orbit.  All other nodes are of 
the trivial single orbit $N$, except the root node having size 256 orbit as shown 
in Section 2.2 
and the depth 1 node having sizes 1 through 8 orbit as  observed in Table~\ref{table:orbits}.  

\begin{figure}
\begin{center}
\includegraphics[width=0.45\textwidth]{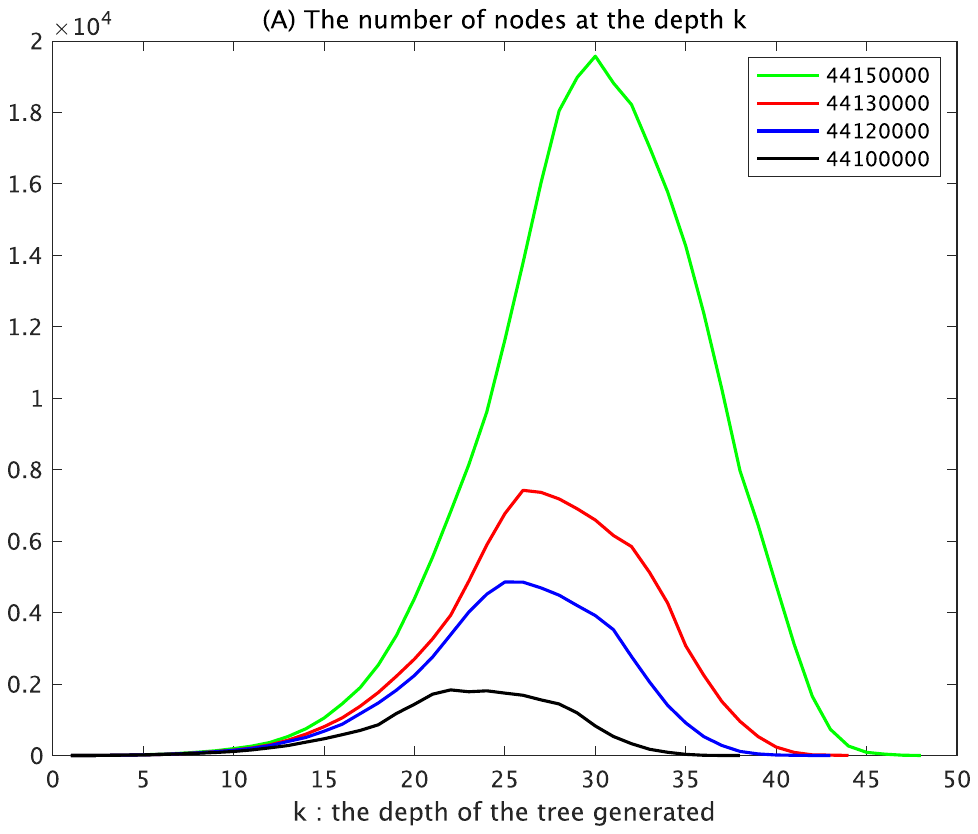}
\includegraphics[width=0.45\textwidth]{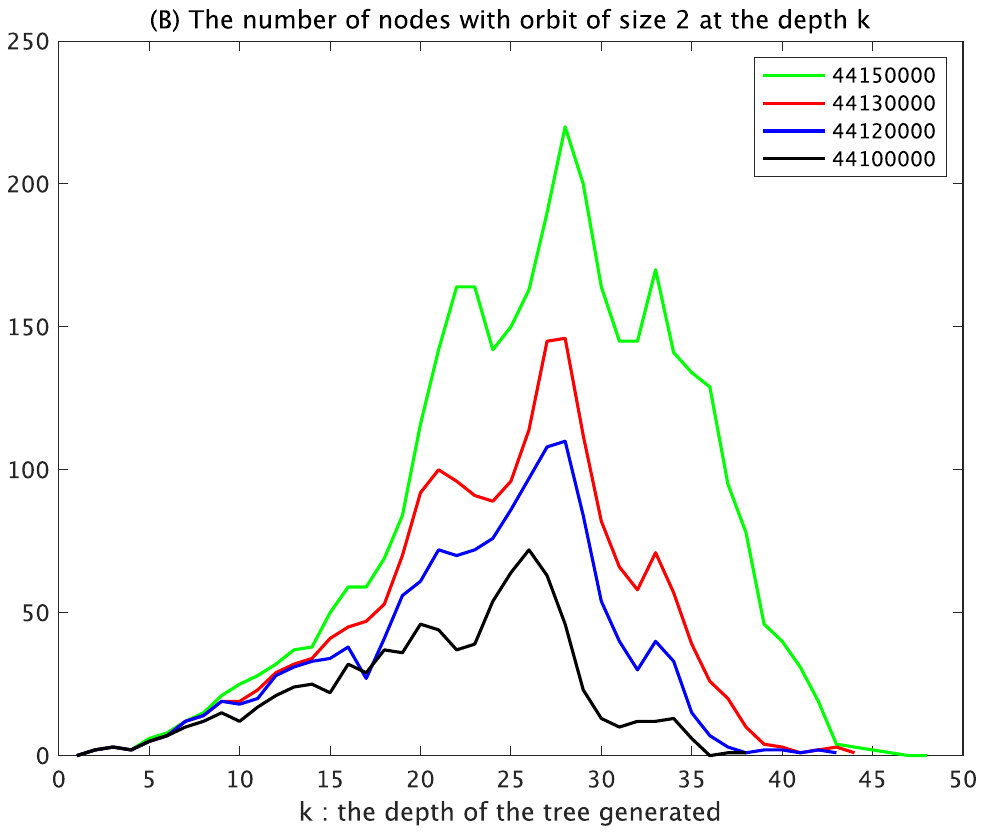}
\vspace{-25mm}
\caption{
(A) 
The number of nodes of the enumeration tree at the depth $k$. 
(B) 
The number of nodes of the enumeration tree with size 2 orbit at the depth $k$. 
}
\end{center}
\end{figure}


\begin{figure}
\begin{center}
\includegraphics[width=0.45\textwidth]{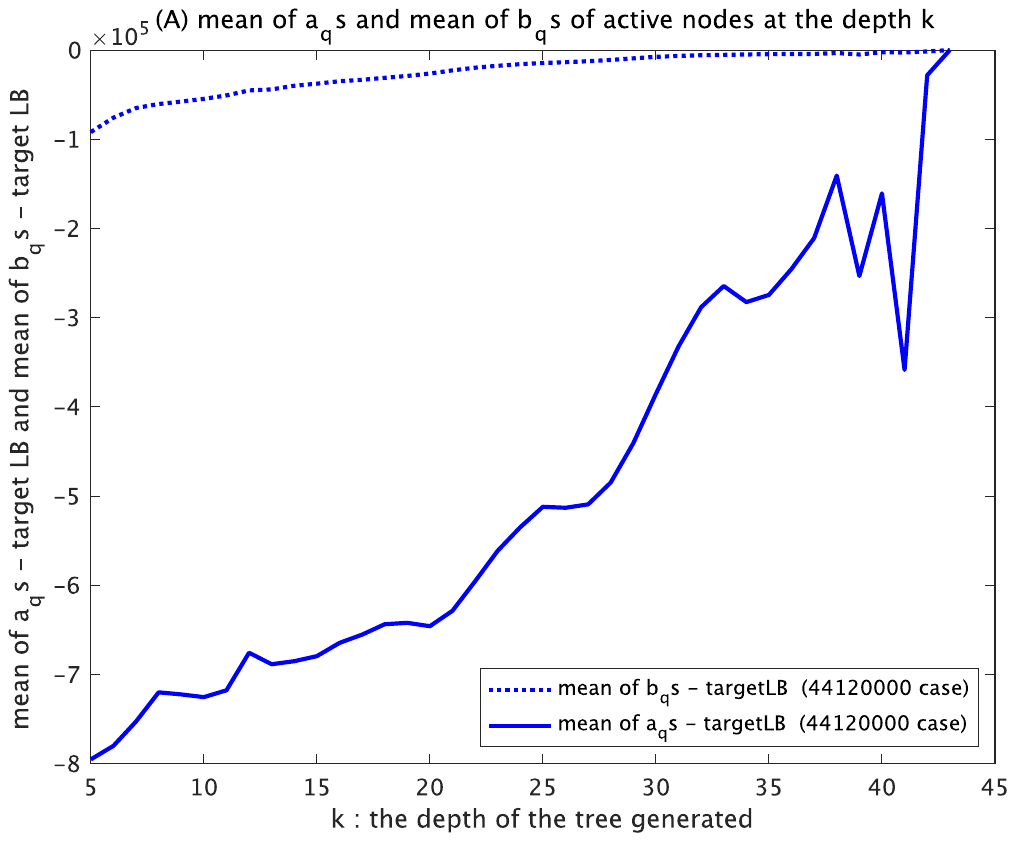} 
\includegraphics[width=0.45\textwidth]{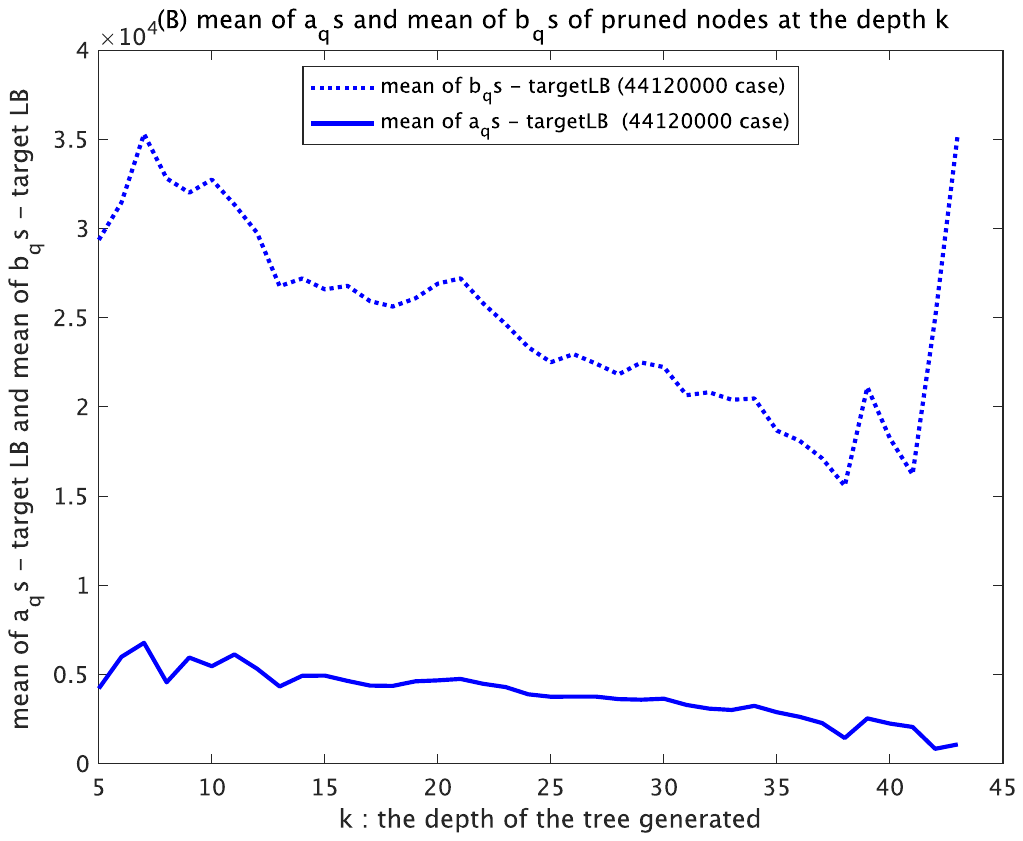} 
\vspace{-25mm}
\caption{The mean of $a_q$ (the blue solid line) and $b_q$ (the blue dotted line) 
when the NB terminated at iteration $q$ 
as $b_q < \hat{\zeta} = 44,120,000$ ({\it i.e.}, active node) --- Case (A) or 
at iteration $q$ as $\hat{\zeta} = 44,120,000 \leq a_q$ ({\it i.e.}, pruned node) --- Case (B).}
\end{center}
\end{figure}


\subsection{Improving the BB method using isomorphism pruning}

\label{section:isomorphism}

Through numerical results  reported in 
Section~4.4, we found that even with the orbital branching,
multiple subproblems appeared in the enumeration tree turned out to be  isomorphic 
(equivalent) to each other. Here, 
two subproblems BQOP$(I_0,I_1)$ and BQOP$(I'_0,I'_1)$ are called {\em isomorphic}  if 
\begin{eqnarray}
& & |I_0| = |I'_0|, \ |I_1| = |I'_1|, \  \sum_{j\in I_1}\sum_{k\in I_1} B_{jk} = \sum_{j\in I'_1}\sum_{k\in I'_1} B_{jk}, \ \mbox{and } \B(I'_0,I'_1) = \P^T\B(I_0,I_1)\P \nonumber \\
& &  \hspace{80mm} \mbox{for some permutation matrix } \P. \label{eq:equivSubproblems}
\end{eqnarray}
The isomorphic subproblems, BQOP$(I_0,I_1)$ and BQOP$(I'_0,I'_1)$, share not only a common optimal value,
  but also a common LB, which is obtained as an optimal value of their 
  Lag-DNN relaxations. 
Therefore, one of them can be pruned even when both of them  are  active. 
This technique is called the isomorphism pruning (see for example, \cite{MARGOT2002} 
and the references therein). 
For a given pair of 
subproblems, BQOP$(I_0,I_1)$ and BQOP$(I'_0,I'_1)$, 
 checking \eqref{eq:equivSubproblems} requires significantly less  CPU time 
 than computing 
their lower bounds. Moreover, various necessary conditions that are easy to implement
can be used for 
verifying \eqref{eq:equivSubproblems}.  
Some of the conditions are, for instance,  
$\sum_{i} [B(I'_0,I'_1)]_{ii} = 
\sum_{i} [B(I_0,I_1)]_{ii}$, $\sum_{i,j} [B(I'_0,I'_1)]_{ij} = 
\sum_{i,j} [B(I_0,I_1)]_{ij}$, 
$\max_{i}[B(I'_0,I'_1)]_{ii} = \max_{i}[B(I_0,I_1)]_{ii}$ and  
$\min_{i,j}[B(I'_0,I'_1)]_{ij} = \min_{i,j}[B(I_0,I_1)]_{ij}$.
By applying those necessary conditions,
the number of the candidates for pairs of 
subproblems BQOP$(I_0,I_1)$ and BQOP$(I'_0,I'_1)$ for which ~\eqref{eq:equivSubproblems} is tested can be considerably reduced, before  verifying  \eqref{eq:equivSubproblems}.

We briefly outline our consistent implementation of the isomorphism 
 pruning with the orbital branching.
Suppose that we have a list of all generated nodes $\{1,\ldots,q\}$ of the depth $1,\ldots,k$ of 
the enumeration tree. The list includes all inactive ({\it i.e.}  already pruned or branched) 
nodes $1,\ldots,p$ located at the depth $1,\ldots,k$ of the enumeration tree for some $p < q$, 
and all active (but not yet branched) nodes $p+1,\ldots,q$ located at the depth $k$.  
We initially employ the orbital branching 
on  active nodes $p+1,\ldots,q$ and generate candidates for  
active nodes to which we will apply the lower bounding procedure, say 
$N_{\rm a} = \{q+1,\ldots,r\}$ for some $r > q$. 
Update the status  of nodes $p+1,\ldots,q$, which have  already undergone branching, from active to inactive.
We then apply the isomorphism pruning 
to the entire node list $\{1,\ldots,r\}$ as follows: 
\begin{description}
\item{Step 1: } Let $N_{\rm d} = \emptyset$ (the set of node to be pruned from the enumeration tree 
by the isomorphism pruning). \vspace{-2mm}
\item{Step 2: } For every $t = q+1,\ldots,r$, if  node $t$ is isomorphic to some node $s < t$,  
then add $t$ to $N_{\rm d}$. \vspace{-2mm}
\item{Step 3: } Update the set $N_{\rm a}$ of active nodes by $N_{\rm a} \backslash N_{\rm d}$.
\vspace{-2mm}
%
\end{description}
Now we are ready to apply the lower bounding procedure to the node set $N_{\rm a}$. 
To simplify the explanation,  
we have assumed that we maintain the list 
of all generated nodes. 
However, only a subset of them is required at Step 2 above. 
For example, let 
$f_1^{\rm min}$ denote the minimum of the numbers of variables fixed to $1$ 
over the node set $N_{\rm a} = \{q+1,\ldots,r\}$. Then, 
any node $s \in \{1,\ldots,q\}$ that has fewer than  $f_1^{\rm min}$ variables fixed to 1 is irrelevant  
at Step 2, 
as well as for future applications of isomorphism pruning, allowing us to eliminate $s$ from the node set 
$\{1,\ldots,q\}$.

Table~4 
shows numerical results on the improved BB method 
in comparison to the original BB method whose numerical results have been 
reported in Section~\ref{section:result}. We observe that the total number of nodes generated in 
the {\bf improved BB} method is less than half of  the one in the original BB method 
in all target LB cases, and that a larger target LB 44,200,000 (1.25\% gap) is newly computed. 

\begin{table}[htp]
\begin{center}
\begin{tabular}{|r|r||r|r|r|r|}
\hline
\multicolumn{2}{|c||}{target LB} & \multicolumn{2}{c|}{no. of nodes} &  \multicolumn{2}{c|}{exec. time (day)} \\
$\hat{\zeta}$ & gap & {\bf improved BB}  & original BB & {\bf improved BB}  & original BB\\ 
\hline
44,100,000  & 1.46\% & {\bf 11,594} & 23,5100 &  {\bf 0.6} &  1.0 \\
44,120,000  & 1.43\% &  {\bf 29,050} & 63,554 &  {\bf 1.2} &  2.5  \\
44,130,000  & 1.41\% &  {\bf 43,904} & 102,310 &  {\bf 1.8} &  4.1  \\
44,150,000  & 1.36\% &  {\bf 109,284} & 277,304 &  {\bf 4.3} &  10.7  \\
44,200,000  & 1.25\% &  {\bf 1,077,353} & 
\multicolumn{1}{c|}{-} &  {\bf 39.2} & \multicolumn{1}{c|}{-} \\ 
\hline
\end{tabular}
\vspace{3mm}
\caption{
Numerical results on the improved BB method in comparison to the original BB method. 
} 
\end{center}
\end{table}

\begin{figure}
\begin{center}
\includegraphics[width=0.45\textwidth]{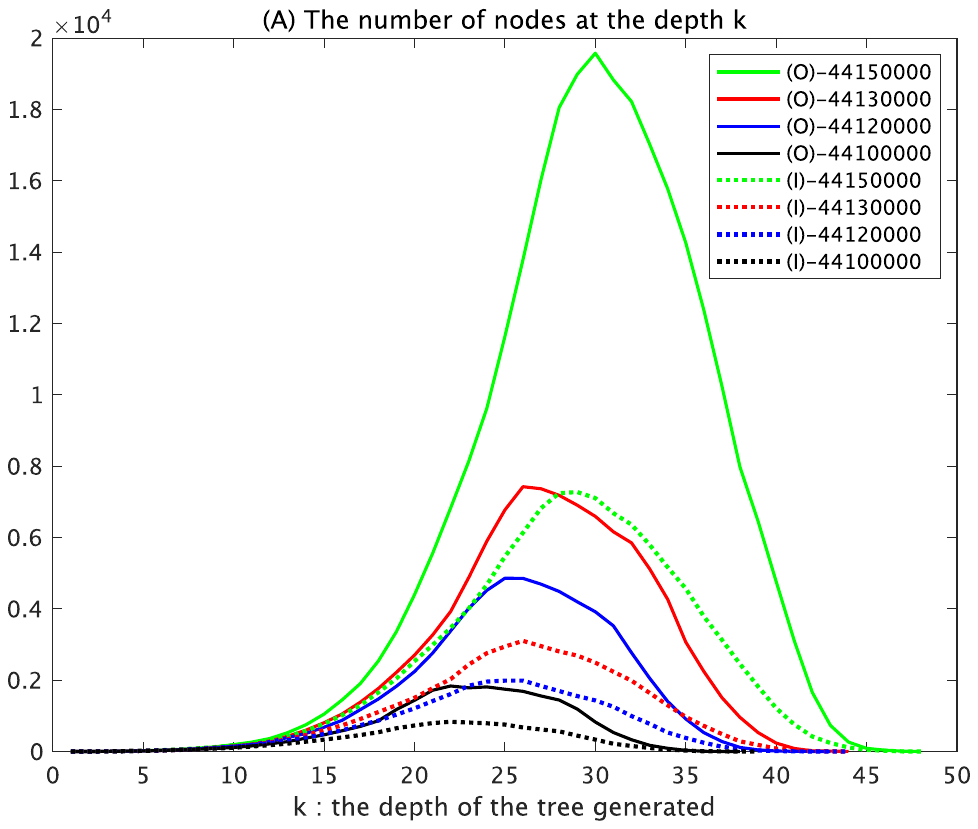}
\includegraphics[width=0.45\textwidth]{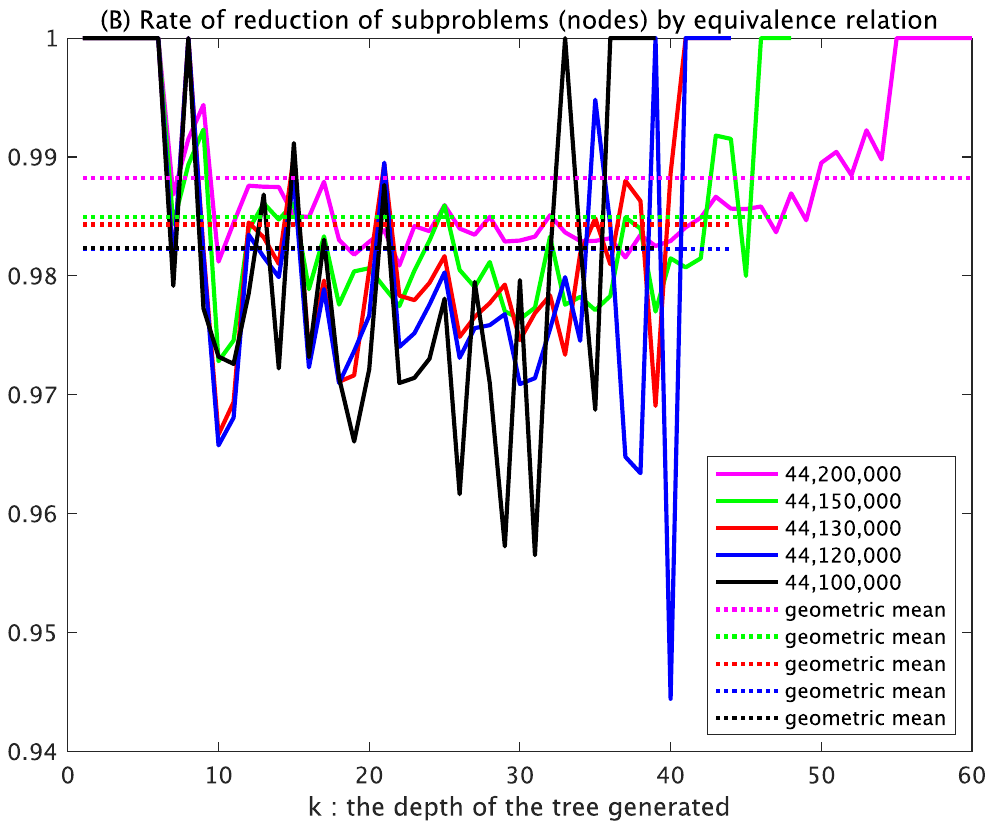}
\vspace{-25mm}
\caption{
(A) Comparison of the numbers of nodes of the enumeration trees in the original BB method (O) and 
 the improved BB method (I) at the depth $k$ of their enumeration trees.  
 (B) Reduction rate of subproblems at the depth $k$ by the equivalence relation.
}
\end{center}
\end{figure}


Figure 4 (A) 
compares the numbers of nodes at the depth $k$ in the original and improved BB methods for 
the target LB = 44,100,000, 44,120,000, 44,130,000 and 44,150,000 cases. We can confirm 
that there exist significant differences between  the numbers of nodes generated by them. 
Figure 4 (B) demonstrates the effectiveness of the 
isomorphism pruning technique  at each depth $k$ 
of the enumeration tree. Let $v_k$ denote the number of active subproblems determined by 
the LB procedure at the depth $k$. Their $2 v_k$ subproblems are generated by 
the orbital branching. By applying the technique, we try to reduce the $2 v_k$ 
subproblems to which the LB procedure is applied at the depth $k+1$. Suppose that 
some $w_k$ nodes are pruned by the technique, 
where $w_k$ could be $0$. Figure 4 (B) shows 
the changes of
$(2v_k-w_k)/(2v_k)$ as $k$ increases $(k=1,\ldots,\ell)$ and their geometric 
mean $r = (\Pi_{k=1}^{\ell}(2v_k-w_k)/(2v_k))^{1/\ell}$, 
where $\ell$ denotes the depth of the 
enumeration tree when the improved BB method terminated. In all target LB cases, 
we see that $r \in [0.98,0.99]$ in Figure 4 (B). 
It can be summarized that 
the technique reduces the number $v_k$ of subproblems generated by the orbital branching 
to $r v_k$ at the depth $k$ on average,  and the modified BB method can reduce  
the total number of 
nodes generated by the original BB method by the factor $r^{\ell}$. 

From the discussions above, we can conclude that the technique proposed 
for pruning equivalent subproblems is indeed effective in accerlerating 
the original BB method. 
We must say, however, that computing an LB with $1.1\%$ gap remains very 
difficult since 
the technique would reduce the number of nodes generated by at most $1/8 \in [0.98,0.99]^{100}$, 
where the improved BB method is assumed to  terminate in the enumeration tree at the  depth $\ell=100$.

\subsection{
Remarks on the difficulty of BQOP~\eqref{eq:BQOP0}
}

\label{section:difficulty}

By incorporating the isomorphism pruning in addition to the orbital branching, 
the improved BB method can entirely avoid the application of the lower bounding and branching procedures to more 
than one isomorphic subproblems of BQOP~\eqref{eq:BQOP0}. 
It should be noted that these two techniques are useful for reducing the number of nodes in 
the enumeration tree that the BB method generates, but not for the lower bounding procedure itself. 
It is apparent that if we could employ a perfectly 
tight lower bounding procedure for any subproblem of BQOP~\eqref{eq:BQOP0}, 
we could immediately determine whether the current UB attains the optimal value. 
However,  this is purely ideal.  
We must acknowledge that even a high-quality lower bounding procedure is  tight only for certain subproblems. Given this reality, 
we can say that solving the 256-dimensional BQOP~\eqref{eq:BQOP0} presents significant challenges,
 in comparison 
to BQOPs of similar size in the benchmark site~\cite{MITTELMANN2023,BIQMAC}. 

When a target LB is selected to be close to the unknown optimal value of BQOP~\eqref{eq:BQOP0}, 
the difficulty in deciding whether to prune each subproblem heavily relies on its optimal value. 
Specifically, 
 if the optimal value is close to the target LB, then the decision becomes more difficult. 
Conversely, if the optimal value is much larger than the target LB, then the decision becomes easier. 
The optimal value of the subproblem is determined by the minimum objective values over 
the feasible solutions of the subproblem. 
Therefore, we expect that the decision is more difficult 
(or easier) as BQOP~\eqref{eq:BQOP0} itself involves more (or fewer) feasible solutions with 
objective values close to the target LB, as the subproblem shares some of them. 
Thus, by sampling,  we investigate the distribution of the objective values 
of feasible solutions of BQOP~\eqref{eq:BQOP0} in 
comparison to other BQOPs of similar size in the benchmark site~\cite{BIQMAC}. 
For the comparison, we apply the following scaling to objective values:
\begin{eqnarray}
\mbox{the scaled objective value} = \frac{\mbox{the objective value} - \mbox{the optimal value of BQOP}}{|\mbox{the optimal value of BQOP}|}. 
\label{eq:scaledObjectiveValue}
\end{eqnarray}
Figure~5 illustrates  
the distribution of the scaled objective values over $1,000,000$ randomly chosen feasible solutions 
for BQOP~\eqref{eq:BQOP0} (the red histogram on the left) and that for 45 BQOPs from the benchmark site~\cite{BIQMAC} 
(the blue histogram
 on the right).  We can observe that the former is clearly 
smaller than the latter. This illustrates why BQOP~\eqref{eq:BQOP0} is notably more challenging 
compared to the benchmark BQOPs. 
 
\begin{figure}
\begin{center}
\includegraphics[width=0.60\textwidth]{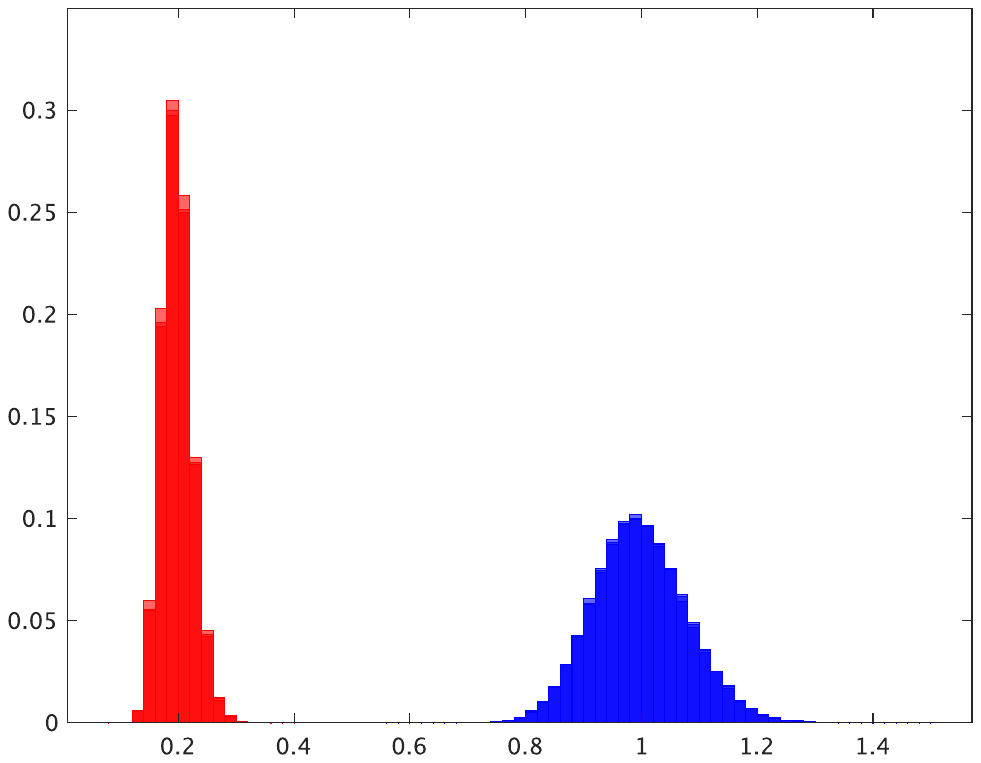} 
\vspace{-40mm}
\caption{
Distribution of the scaled objective values defined 
by~\eqref{eq:scaledObjectiveValue}. The horizontal axis represents the scaled 
objective value and the vertical axis the probability. 
The red histogram on the left depicts the distribution of 
 the scaled objective values of  1,000,000 randomly chosen feasible solutions of BQOP~\eqref{eq:BQOP0}. 
 The blue 
 histogram on the right  depicts that of the 45 BQOPs (bqp250-1,$\ldots$,bqp250-10, be200.3.1,$\ldots$,be200.3.10,
be200.8.1,$\ldots$,be200.8.10,be250.1$\ldots$,be250.10,gka1e,$\ldots$,gka5e) with dimensions  $200$ or $250$, where 1,000,000 feasible solutions are randomly sampled in each BQOP.
}
\end{center}
\end{figure}


\section{Concluding remarks}

\label{section:concludingRemarks} 

We have investigated the $256$-dimensional BQOP with a single 
cardinality constraint, BQOP~\eqref{eq:BQOP0}, which is converted from the 
largest unsolved QAP instance tai256c. 
The converted BQOP with dimension $256$ is much simpler than the original QAP 
tai256c involving $256 \times 256 = 65536$ binary variables, 
and its dimension $256$ is not so large.
While one might expect the converted BQOP to be notably easier to solve compared to the original QAP tai256c,
our findings indicate that it still presents a significant challenge. 
The challenge
primarily stems from the symmetry property~\eqref{eq:symmetry3} exhibited in the coefficient matrix $\B$, which is inherited from tai256c.
For future development toward solving the BQOP, 
we need \vspace{-2mm}
\begin{itemize}
\item an efficient and much stronger lower bounding procedure than the DNN relaxation, 
\vspace{-2mm}
\item additional  techniques  to enhance the exploitation of the symmetry property~\eqref{eq:symmetry3}, and  
\vspace{-2mm}
\item more powerful computer systems.
\end{itemize}

While we have focused on BQOP~\eqref{eq:BQOP0}  
converted from the QAP tai256c in this paper, it is straightforward to adapt the discussion of the paper  
to  general BQOPs with a single cardinality 
constraint and  general QUBOs which satisfy the symmetry 
property~\eqref{eq:symmetry3}. 

\subsection*{Acknowledgements}
The authors are grateful to Professor Koji Nakano for providing them with 
numerical results on 
DABS (Diverse Adaptive Bulk Search, a genetic algorithm-based search algorithm 
for solving QUBO \cite{NAKANO2023}), which have been included in Section~4.
They are also grateful to Tobias Achterberg for sharing the internal information of Gurobi and to Christopher Hojny for discussions about symmetry.
The work of Sunyoung Kim was supported by  NRF 2021-R1A2C1003810.
The work for this article has been partially conducted within the Research Campus MODAL funded by the German Federal Ministry of Education and Research (BMBF grant numbers
05M14ZAM, 05M20ZBM).


\end{document}